\documentclass[12pt,oneside]{amsproc}
\usepackage[T2A]{fontenc}
\usepackage[cp1251]{inputenc}
\usepackage[russian]{babel}
\usepackage{amssymb}
\usepackage{amsthm}

\newtheorem{thm}{Теорема}
\newtheorem{deff}{Определение}
\newtheorem{statement}{Утверждение}
\newtheorem{cons}{Следствие}
\newtheorem{lemma}{Лемма}
\newtheorem*{example}{Пример}

\renewcommand{\span}{\operatorname{span}}
\renewcommand{\Im}{\operatorname{Im}}
\renewcommand{\Re}{\operatorname{Re}}

\newcommand{\opi}{\operatorname{I}}
\newcommand{\opv}{\operatorname{V}}
\newcommand{\opt}{\operatorname{t}}
\newcommand{\itL}{\mathit{L}}

\title{Диссипативные краевые условия для обыкновенных дифференциальных
операторов}
\author{Е.~А.~Ширяев}
\thanks{Работа выполнена при поддержке грантов РФФИ 04-01-00712 и НШ-1927.2003.1.}
\date{}
\begin{document}
\begin{flushleft}
УДК 517.984
\end{flushleft}
\maketitle

В работе \cite{Birkhoff} Дж.~Д.~Биркгоф развил асимптотические
методы для исследования обыкновенных дифференциальных операторов
высокого порядка, порождённых выражением

\begin{equation}
\label{begeq} l(y)=(-1)^m
y^{(m)}+(-1)^{m-2}[p_2(x)y]^{(m-2)}+\cdots+p_m(x)y,
\end{equation}
где функции $p_k(x)$ бесконечно дифференцируемы на $[a,b]$, и
краевыми условиями
\begin{equation}
\label{begbc} U_j(y)=\sum \limits_{k=0}^{m-1}
a_{j\,k}y^{(k)}(a)+b_{j\,k}y^{(k)}(b)=0,\; j=1,\ldots ,m.
\end{equation}

В частности, он выделил важный класс краевых условий, которые
назвал регулярными. Оператор, порождённый выражением \eqref{begeq}
и регулярными краевыми условиями \eqref{begbc}, здесь мы также
условимся называть регулярным. Отметим, что  в определении
регулярности участвуют только коэффициенты в линейных формах
\eqref{begbc}, и не участвуют коэффициенты $p_j(x)$
дифференциального выражения.

Важный результат, установленный Биркгофом, состоял в оценке
резольвенты регулярного дифференциального оператора. Оценка
получалась, по существу, такой же, как для самосопряжённых
краевых условий, т.е. условий, при которых оператор $\itL_0$,
порождённый выражением ${l_0(y)=(-i)^my^{(m)}}$ самосопряжён.
Задача о том, являются ли самосопряжённые условия регулярными,
оказалась непростой. Она была положительно решена С.~Салафом
\cite{Sal} для чётных $m$ и для произвольного порядка
А.~М.~Минкиным \cite{Minkin_MZ}.

В этой заметке расматривается оператор $\itL_0$, порождённый
дифференциальным выражением ${l_0(y)=(-i)^my^{(m)}}$ c областью
опеределения ${D(\itL_0)=\{y \in W_2^m[0,1], \, \mathit{U}_j(y)=0,
\, j=1,\ldots,m\}}$. Здесь через $W_2^m$ обозначено пространство
Соболева.

Краевые условия \eqref{begbc} назовём диссипативными, если
оператор $\itL_0$ является диссипативным, т.е. $\Im (\itL_0y,\,y)
\geqslant 0$.

Целью работы является описание всех диссипативных краевых условий
и выяснение вопроса, являются ли такие условия регулярными.
Отметим следующий факт: если мы докажем, что диссипативный
оператор $\itL_0$ регулярен, то произвольный оператор $\it L$,
порождённый дифференциальным выражением \eqref{begeq} и теми же
краевыми условиями, будет регулярен.

\section{Описание диссипативных краевых условий.}
В этом параграфе приведены необходимые и достаточные условия
диссипативности краевых условий.

Для краткости записи введём векторы-строки:
$$\widehat{y}_0:=(y(0), \ldots, y^{(m-1)}(0));$$
$$\widehat{y}_1:=(y(1), \ldots ,y^{(m-1)}(1));$$
$$\widehat{y}:=(\widehat{y}_0,\; \widehat{y}_1)
=\bigl(y(0),\ldots, y^{(m-1)}(0),\;y(1), \ldots,
y^{(m-1)}(1)\bigr).$$

Через $\widehat{y}_0^{\,*}$, $\widehat{y}_1 ^{\,*}$,
$\widehat{y}^{\,*}$, $A^*$ будем обозначать соответствующие
векторы-столбцы и матрицы с комплексно сопряжёнными компонентами,
а символ \glqq${\opt}$\grqq~на месте верхнего индекса будет
означать транспонирование вектора или матрицы без комплексного
сопряжения компонент. Для обозначения матриц и векторов,
компоненты которых комплексно сопряжены исходным, будем ставить
сверху черту. То есть, $A^*=\overline{A} ^{\,\opt}$

Заметим, что краевые условия в новых обозначениях можно
переписать в виде
${A\widehat{y}_0^{\,\,\opt}+B\widehat{y}_1^{\,\,\opt}=0}$, где $A$
и $B$---матрицы порядка $m$.

\begin{lemma} \label{P} Имеем
$$2\Im{(\mathit{L}_0y,y)}=
\left\{ \begin{array}{rcl}{i(-1)^{n+1}\bigl[ \,\widehat{y}_1
\mathit{J}\widehat{y}_1^{\,*}-
\widehat{y}_0 \mathit{J}\widehat{y}_0^{\,*} \,\bigr],\,\,m=2n}\vspace{1.5ex}\\
{(-1)^n\bigl[\,\widehat{y}_1\mathit{K}\widehat{y}_1^{\,*}-
\widehat{y}_0\mathit{K}{\widehat{y}_0}^{\,*}\,\bigr], \,\,m=2n-1,}
\end{array}
\right.$$ где матрицы $\mathit{J}$ и $\mathit{K}$ выглядят
следующим образом:
\smallskip

$$ \mathit{J}=
\begin{pmatrix}
0 & 0 & \ldots & 0 & -1\\
0&&& 1 & 0\\ \vdots && \cdot && \vdots \\
0 & -1 &&& 0 \\ 1 & 0 & \ldots && 0
\end{pmatrix} \qquad
\mathit{K}=
\begin{pmatrix}
0 & 0 & \ldots & 0 & 1\\
0&&& -1 & 0\\ \vdots && \cdot && \vdots \\
0 & -1 &&& 0 \\ 1 & 0 & \ldots && 0
\end{pmatrix}
$$
\end{lemma}
\begin{proof}

$$\Im{(\mathit{L}_0y,y)}=\frac{1}{2i}\bigl( (\mathit{L}_0y,y)-(y,\mathit{L}_0y)\bigr)$$
Интегрируя $(\mathit{L}_0y,y)$ $m$ раз по частям, получим
$$2\Im{(\mathit{L}_0y,y)}=(-i)\bigl[(-i)^{m}\bigl(\left.y^{(m-1)}\overline{y} -
\ldots +(-1)^{m-1}y\overline{y^{(m-1)}}\bigr) \right|_0^1\bigr]=$$
$$=
\left\{ \begin{array}{rcl}{i(-1)^{n+1}\bigl[\,\widehat{y}_1
\mathit{J}\widehat{y}_1^{\,*}-
\widehat{y}_0 \mathit{J}\widehat{y}_0^{\,*}\,\bigr],\,\,m=2n}\vspace{1ex}\\
{(-1)^n\bigl[\,\widehat{y}_1\mathit{K}\widehat{y}_1^{\,*}-
\widehat{y}_0\mathit{K}{\widehat{y}_0}^{\,*}\,\bigr], \,\,m=2n-1,}
\end{array}
\right.$$

\end{proof}

\begin{lemma} \label{lm4}
Пусть $m=2n$ и краевые условия заданы в форме
${A\widehat{y}_0^{\,\,\opt}+B\widehat{y}_1^{\,\,\opt}=0}$, где
$A$ и $B$---матрицы порядка $2n$. Если $\forall y \in
\mathit{D}(\mathit{L})$
$$
\widehat{y}_1(-i(-1)^n\mathit{J})\widehat{y}_1^{\,*} -
\widehat{y}_0(-i(-1)^n\mathit{J})\widehat{y}_0^{\,*} \geqslant 0,
$$
то

\begin{equation}
\label{disscond}
[\overline{A}(-i(-1)^n\mathit{J})A^{\opt}-\overline{B}(-i(-1)^n\mathit{J})B^{\opt}]
\leqslant 0.
\end{equation}

\end{lemma}

\begin{proof}[Доказательство]

Обозначим $\widehat{U}_j=(a_{0j},\; a_{1j},\ldots, a_{2n-1\,j},\;
b_{0j},\; b_{1j},\ldots, b_{2n-1\,j}), \; j=1,\ldots,2n$ строки,
состоящие из коэффициентов краевых условий \eqref{begbc}. Тогда
равенство $U_j(y)=0$ можно переписать в виде

\begin{equation}
\label{condition}  \langle\overline{\widehat{U}}_j,\;
\widehat{y}^{\,\,\opt}\rangle_{4n}=0
\end{equation}

Последнее выражение --- скалярное произведение в $4n$ - мерном
пространстве над $\mathbb{C}$. Строки $\overline{\widehat{U}}_j$
образуют линейное подпространство размерности $2n$ (здесь и далее
будем полагать, что краевые условия линейно независимы). Обозначим
через $\mathfrak{N}$ подпространство размерности $2n$ состоящее
из решений системы уравнений \eqref{condition}, а через
$\mathfrak{L}$
--- подпространство $\span(\{\overline{\widehat{U}}_j\}_{j=1}^{2n})$.
Тогда $\mathfrak{N}$ ортогонально $\mathfrak{L}$. Положим
$$ \mathit{M}:=\left( \begin{array}{cc} i(-1)^n\mathit{J}&0 \\ 0 &
-i(-1)^n\mathit{J}
\end{array} \right), $$
тогда
\begin{equation}
\label{interpretation} \widehat{y} \, \mathit{M} \, \widehat{y}^*
\geqslant 0 \quad \forall\, \widehat{y} \in \mathfrak{N}.
\end{equation}

\medskip
У матрицы $\mathit{M}$  два собственных значения: $1$ и $(-1)$.
Каждое из них имеет кратность $2n$. Пусть
$\{\widehat{x}_j^{\,-}\}_{j=1}^{2n}, \; \;
\{\widehat{x}_j^{\,+}\}_{j=1}^{2n}$
 --- базис, в котором $\mathit{M}$
диагонализируется: сперва на диагонали стоят
\mbox{\glqq$-1$\grqq}, потом \mbox{\glqq$1$\grqq}. Перейдём в
$\mathbb{C}^{4n}$ посредством ортогонального преобразования к
базису, состоящему из векторов
$\{\widehat{x}_1^{\,-},\,\widehat{x}_2^{\,-},\ldots,\widehat{x}_{2n}^{\,-},\,
\widehat{x}_1^{\,+},\,\widehat{x}_2^{\,+},\ldots,\widehat{x}_{2n}^{\,+}\}$.
Обозначим
$\mathfrak{N}_-=\span(\{\widehat{x}_j^{\,-}\}_{j=1}^{2n})$ и
$\mathfrak{N}_+=\span(\{\widehat{x}_j^{\,+}\}_{j=1}^{2n})$.
Отметим, что $\forall \, \widehat{x} \in \mathfrak{N}_- \;\;
\widehat{x}M\widehat{x}^{\,*}\leqslant 0$ и $\forall \,
\widehat{x} \in \mathfrak{N}_{\,+} \;\;
\widehat{x}M\widehat{x}^{\,*}\geqslant 0$.

Рассмотрим операторы $\mathit{S}_-$ и $\mathit{S}_+$, которые
проецируют $\mathfrak{N}$ на $\mathfrak{N}_-$ и $\mathfrak{N}_+$
соответственно. $\forall \, \widehat{x} \in \mathfrak{N}_+$
существует не более одного $\widehat{y} \in \mathfrak{N}: \;
\mathit{S}_+(\widehat{y})=\widehat{x}$. В противном случае, если
найдутся такие $\widehat{y}_1$ и $\widehat{y}_2$, ${(\widehat{y}_1
- \widehat{y}_2)\mathit{M}(\widehat{y}_1 - \widehat{y}_2)^*<0}$.
Это означает, что у оператора $\mathit{S}_+$ есть обратный (также
линейный). В силу этого подпространству $\mathfrak{N}$
принадлежат те и только те векторы, координаты которых в базисе
$\{\widehat{x}_1^{\,-},\,\widehat{x}_2^{\,-},\ldots,\widehat{x}_{2n}^{\,-},\,
\widehat{x}_1^{\,+},\,\widehat{x}_2^{\,+},\ldots,\widehat{x}_{2n}^{\,+}\}$
состоят  из последовательно выписанных координат вектора
$(\mathit{S}_-\mathit{S}_+^{-1})(\widehat{x})$ в базисе
$\{\widehat{x}_j^{\,-}\}_{j=1}^{2n}$ пространства
$\mathfrak{N}_-$ и координат $\widehat{x}$ в базисе
$\{\widehat{x}_j^{\,+}\}_{j=1}^{2n}$ пространства $\mathfrak{N}_+$
($\forall \widehat{x}\in \mathfrak{N}_+$). Будем записывать это
так:
$$
\mathfrak{N}= \big
\{\big((\mathit{S}_-\mathit{S}_+^{-1})(\widehat{x})\,,\,\widehat{x}\,\big)\mid
\;\widehat{x} \in \mathfrak{N}_+\big \}.
$$

Теперь дадим аналогичное описание для
$\mathfrak{L}=\mathfrak{N}^{\perp}$. А именно, представим
оператор ${\mathit{T}: \mathfrak{N}_- \rightarrow
\mathfrak{N}_+}$ такой, что
$\mathfrak{L}=\mathfrak{N}^{\perp}=\big
\{\big(\,\widehat{y},\,\mathit{T}(\widehat{y})\big)
\mid\;\widehat{y} \in \mathfrak{N}_-\big \}$. Зафиксируем в
подпространствах $\mathfrak{N}_-$ и $\mathfrak{N}_+$ базисы
$\{\widehat{x}_j^{\,-}\}_{j=1}^{2n}$ и
$\{\widehat{x}_j^{\,+}\}_{j=1}^{2n}$ соответственно. Пусть
$S$---матрица оператора ${\mathit{S}_-\mathit{S}_+^{-1}:
\mathfrak{N}_+ \rightarrow \mathfrak{N}_-}$ в этих базисах. Тогда
рассмотрим линейный оператор из $\mathfrak{N}_-$ в
$\mathfrak{N}_+$, матрица которого при фиксированных базисах
$\{\widehat{x}_j^{\,-}\}_{j=1}^{2n}$ и
$\{\widehat{x}_j^{\,+}\}_{j=1}^{2n}$ равна ${(-S^*)}$. Это и есть
искомый оператор $\mathit{T}$.

Из \eqref{interpretation} следует, что ${\forall \, \widehat{x}
\in \mathfrak{N}_+}\; \langle
(\mathit{S}_-\mathit{S}_+^{-1})(\widehat{x}) ,
\,(\mathit{S}_-\mathit{S}_+^{-1})(\widehat{x}) \rangle_{2n}
\leqslant \langle \widehat{x} , \, \widehat{x}\rangle_{2n}$, а
значит, ${\forall \, \widehat{x} \in \mathfrak{N}_-}\; \langle
\mathit{T}(\widehat{x}), \, \mathit{T}(\widehat{x}) \rangle_{2n}
\leqslant \langle \widehat{x} , \, \widehat{x}\rangle_{2n}$.
Следовательно, $\forall \, \widehat{y} \in \mathfrak{N}^{\perp}
\quad \widehat{y} \,\mathit{M}\, \widehat{y}^{\,*} \leqslant 0$.

Поскольку $\mathfrak{N}^{\perp}=\mathfrak{L}$---подпространство,
являющееся линейной оболочкой
$\{\overline{\widehat{U}}_j\}_{j=1}^{2n}$, то
$$[\overline{A}(-i(-1)^n\mathit{J})A^{\opt}-\overline{B}(-i(-1)^n\mathit{J})B^{\opt}] \leqslant 0$$
Лемма \ref{lm4} доказана.

\end{proof}

Из доказательства леммы \ref{lm4} ясно, что если $A$ и $B$
удовлетворяют \eqref{disscond}, то они задают диссипативные
условия.

В нечётном случае диссипативные условия описываются так же.
Причина этого в том, что и в нечётном случае у матрицы
$\mathit{M}$ будет одинаковое количество \mbox{\glqq$1$\grqq} и
\mbox{\glqq$-1$\grqq}~в диагонализованном виде. Таким образом,
доказана теорема, дающая представление об алгебраической структуре
диссипативных расширениях оператора $\mathit{L}_0$.

\begin{thm}

Матрицы $A$ и $B$ задают диссипативные краевые условия тогда и
только тогда, когда
$$
[\overline{A}(-i(-1)^n\mathit{J})A^{\opt}-\overline{B}(-i(-1)^n\mathit{J})B^{\opt}]
\leqslant 0
$$
в чётном случае,

$$
[\overline{A}(-1)^n\mathit{K}A^{\opt}-\overline{B}(-1)^n\mathit{K}B^{\opt}]
\leqslant 0
$$
в нечётном случае.
\end{thm}

\bigskip
Предложим теперь ещё одно описание диссипативных краевых условий
для дифференциального выражения ${l_0(y)=(-i)^my^{(m)}}$. Этот
результат является развитием работы Ф.С.Рофе-Бекетова и
А.М.Холькина \cite{RH}.

Начнём с чётного случая: ${m=2n}$. Для удобного выражения формы
${\Im{(\mathit{L}_0y,y)}}$ нам потребуются квазипроизводные,
отвечающие дифференциальному выражению ${l_0(y)}$. Определим их
следующим образом:
\begin{flushleft}
\quad\quad$y^{[k]}=y^{(k)}, \; (k=0,1,\ldots,n;\;y^{(0)}=y),$\\
\quad\quad$y^{[2n-k]}=-\frac{d}{dt}y^{[2n-k-1]}, \;
(k=0,1,\ldots,n-1).$
\end{flushleft}
Введём обозначения для векторов младших и старших производных:
\begin{flushleft}
\quad\quad$y^{\land}=\bigl(y(0),\; y'(0),\ldots,
y^{(n-1)}(0),\;y(1),\;
y'(1),\ldots, y^{(n-1)}(1) \bigr)^{\opt},$\\
\quad\quad$y^{\lor}=\bigl(y^{[2n-1]}(0),\; y^{[2n-2]}(0),\ldots,
y^{[n]}(0),\;-y^{[2n-1]}(1),\; -y^{[2n-2]}(1),\ldots, -y^{[n]}(1)
\bigr)^{\opt}.$
\end{flushleft}
С помощью этих обозначений ${\Im{(\mathit{L}_0y,y)}}$ можно
записать так:
$$\Im{(\mathit{L}_0y,y)}=\frac{1}{2i}\bigl((\mathit{L}_0y,y)-
(y,\mathit{L}_0y)\bigr)=\frac{1}{2i}\bigl(\langle
y^{\lor},y^{\land}\rangle_{2n}- \langle
y^{\land},y^{\lor}\rangle_{2n}\bigr),$$
где ${\langle \cdot\, ,
\cdot\rangle_{2n}}$ по-прежнему означает скалярное произведение в
$\mathbb{C}^{2n}$.

\begin{lemma} \label{RH_even1}
Если краевые условия задаются в виде
$$(\opv-\opi)y^{\lor}+i(\opv+\opi)y^{\land}=0,$$
где $\opi$---единичный, а $\opv$---нерастягивающий линейный
оператор в $\mathbb{C}^{2n}$, то есть
$$\| \opv{y} \|\leqslant
\|y\|\; \;  \forall y \in \mathbb{C}^{2n},$$

то оператор, порождённый выражением $l_0(y)$ и этими краевыми
условиями, диссипативен.
\end{lemma}
\begin{proof}
Перепишем краевые условия в виде
$$\opv(y^{\lor}+iy^{\land})=y^{\lor}-iy^{\land},$$
тогда
$$\|y^{\lor}+iy^{\land} \|^2\geqslant \|y^{\lor}-iy^{\land}\|^2.$$
Отсюда
$$2i\bigl[ \langle y^{\land},y^{\lor}\rangle_{2n}-\langle y^{\lor},y^{\land}\rangle_{2n} \bigr]\geqslant 0,$$
что совпадает с определением диссипативности оператора.
\end{proof}

\begin{lemma} \label{RH_even2}
Пусть оператор $\mathit{L}_0$ порождён дифференциальным выражением
$l_0(y)$, а его область определения ${D(\mathit{L}_0)\subset
W_2^{2n}(0,1)}$ задаётся $2n$ линейно независимыми краевыми
условиями. Если $\mathit{L}_0$ диссипативен, то существует такой
нерастягивающий линейный оператор $\opv$, действующий в
$\mathbb{C}^{2n}$, что краевые условия эквивалентны следующим:
$$(\opv-\opi)y^{\lor}+i(\opv+\opi)y^{\land}=0,$$
где $\opi$ по-прежнему единичный оператор в ${\mathbb{C}^{2n}}$.
\end{lemma}
\begin{proof}
Поскольку оператор $\mathit{L}_0$ диссипативен, то
$$2i\bigl[ \langle y^{\land},y^{\lor}\rangle_{(2n)}-\langle y^{\lor},y^{\land}\rangle_{2n} \bigr]\geqslant 0,$$
что равносильно
$$\|y^{\lor}+iy^{\land} \|^2\geqslant \|y^{\lor}-iy^{\land}\|^2.$$
Определим теперь отображение $\opv$.

\noindent Область определения $D(\opv)$ зададим так:
$$D(\opv)=\{z\in \mathbb{C}^{2n}|\, z=y^{\lor}+iy^{\land}, \; y\in D(\mathit{L}_0)\}.$$
Из условия леммы следует что $D(\opv)$---подпространство в
$\mathbb{C}^{2n}$. Положим
$$\opv(y^{\lor}+iy^{\land})=y^{\lor}-iy^{\land}, \; y\in D(\mathit{L}_0).$$
Отображение $\opv: \, \mathbb{C}^{2n}\longrightarrow
\mathbb{C}^{2n}$ является искомым линейным оператором.

Докажем, что заданное так отображение $\opv$ однозначно.

Предположим противное: пусть отображение $\opv$ неоднозначно, то
есть найдутся такие ${q,\: y\in D(\mathit{L}_0), \; q \neq y}$,
что ${y^{\lor}+iy^{\land}= q^{\lor}+iq^{\land}}$, и
${y^{\lor}-iy^{\land}\neq q^{\lor}-iq^{\land}}$. Тогда, так как
${(y-q)\in D(\mathit{L}_0)}$, то
$$0= \|(y-q)^{\lor}+i(y-q)^{\land} \| \geqslant
\|(y-q)^{\lor}-i(y-q)^{\land}\| >0.$$ Пришли к противоречию с
предположением о неоднозначности отображения $\opv$.

Таким образом, представлен линейный нерастягивающий оператор
$\opv$, действующий в $\mathbb{C}^{2n}$, такой, что краевые
условия
$$(\opv-\opi)y^{\lor}+i(\opv+\opi)y^{\land}=0$$
задают ${D(\itL _0)}$. Доказательство леммы \ref{RH_even2}
закончено.
\end{proof}

Для нечётного случая (${m=2n-1}$) определение квазипроизводных
$y^{[k]}$, отвечающих ${l_0(y)}$, задаётся формулами:

\begin{flushleft}
\quad\quad$y^{[k]}=y^{(k)}, \; (k=0,1,\ldots,n-2;\;y^{(0)}=y),$\\
\quad\quad$y^{[n-1]}=-iy^{(n-1)},$\\
\quad\quad$y^{[2n-k-1]}=-\frac{d}{dt}y^{[2n-k-2]}, \;
(k=0,1,\ldots,n-1).$
\end{flushleft}

Обозначим:
\begin{flushleft}
\quad\quad$y^{\land}=\bigl(y^{(n-1)}(0)+y^{(n-1)}(1),\; y(0),\;
y'(0),\ldots, y^{(n-2)}(0),\;y(1),\; y'(1),\ldots, y^{(n-2)}(1)
\bigr)^{\opt},$\\
\quad\quad$y^{\lor}=\bigl(iy^{(n-1)}(1)-iy^{(n-1)}(0),\;y^{[2n-2]}(0),\ldots,
y^{[n]}(0),\;-y^{[2n-2]}(1),\ldots, -y^{[n]}(1) \bigr)^{\opt}.$
\end{flushleft}

Теперь ${\Im(\mathit{L}_0y,y)}$ снова можно записать так:
$$\Im{(\mathit{L}_0y,y)}=\frac{1}{2i}\bigl(\langle y^{\lor},y^{\land}\rangle_{2n-1}-
\langle y^{\land},y^{\lor}\rangle_{2n-1}\bigr),$$ где ${\langle
\cdot\, , \cdot\rangle_{2n-1}}$ означает скалярное произведение в
$\mathbb{C}^{2n-1}$.

Леммы \ref{RH_even1} и \ref{RH_even2} верны в нечётном случае. Их
совокупность доказывает следующую теорему.

\begin{thm}
Оператор, порождённый ${l_0(y)=(-i)^my^{(m)}}$ и краевыми
условиями \eqref{begbc}, является  диссипативным тогда и только
тогда, когда найдётся такой нерастягивающий оператор $\opv$ в
$\mathbb{C}^{m}$, что можно переписать эти условия в
эквивалентной форме:
$$(\opv-\opi)y^{\lor}+i(\opv+\opi)y^{\land}=0.$$
\end{thm}

\section{Определение регулярных краевых условий}
Напомним определение регулярности краевых условий (см.
\cite{Birkhoff}, \cite{Naimark}). Предварительно краевые условия
нужно нормировать. Число $k_j$ ($j=1,\ldots ,m$) назовём порядком
краевого условия ${U_j(y)=0}$, если это краевое условие содержит
${y^{(k_j)}(0)}$ или ${y^{(k_j)}(1)}$, но не содержит
${y^{(\nu)}(0)}$ и ${y^{(\nu)}(1)}$ при $\nu>k_j$. Рассмотрим
краевые условия порядка ${m-1}$, если такие имеются. Заменяя их,
если надо, линейными комбинациями, можно добиться того, чтобы
число линейно независимых краевых условий порядка ${m-1}$ было
$\leqslant 2$. Остальные краевые условия имеют порядок
${\leqslant m-2}$; применяя к ним тот же приём, сведём их число к
минимуму и т.д.

Описанные операции называются {\it нормировкой краевых условий}, а
полученные в результате краевые условия называются {\it
нормированными}. Из способа их построения следует, что
нормированные краевые условия должны иметь вид:

\begin{equation}
\label{normbc} U_j(y)\equiv U_{j\,0}(y)+U_{j1}(y)=0,
\end{equation}
 где
$$U_{j\,0}(y)=\alpha_j y^{(k_j)}(0)+\sum
\limits_{s=0}^{k_j-1}\alpha_{j\,s}y^{(s)}(0), \quad
U_{j1}(y)=\beta_j y^{(k_j)}(1)+\sum
\limits_{s=0}^{k_j-1}\beta_{j\,s}y^{(s)}(1),$$
$$m-1 \geqslant k_1 \geqslant k_2\geqslant \ldots \geqslant k_m \geqslant 0,\; k_{j} > k_{j+2},$$
причём для каждого значения индекса $j$ хотя бы одно из чисел
${\alpha_j,\:\beta_j}$ отлично от $0$.

Отметим, что имеется эквивалентное определение нормированных
краевых условий (см. \cite{shef}).
\begin{deff}
Будем говорить, что краевые условия $V_j(y)=0$, $j=1,\ldots , m$,
эквивалентные \eqref{begbc}, реализуют условия \eqref{begbc} в
минимальной форме (или приводят их к минимальной форме), если
сумма порядков условий $V_j(y)$ минимальна.
\end{deff}

\begin{statement}
Краевые условия нормированы тогда и только тогда, когда они
приведены к минимальной форме.
\end{statement}
\begin{proof}
Очевидно, из минимальности формы следует нормированность.
Действительно, пусть, например, условия $\{{U}_j(y)=0\}_{j=1}^m$
приведены к минимальной форме, но не нормированы. Это означает,
что найдутся два или более условия некоторого порядка $k$ такие,
что есть линейная комбинация строк длины два
$\sum\limits_{j}c_j(a_{jk},\;b_{jk})=(0,\; 0)$. Здесь $j$
принимает значения, равные номерам условий порядка $k$. Но тогда
у одного из этих условий порядок можно понизить, а это приводит к
противоречию с предположением о минимальности формы
$\{{U}_j(y)=0\}_{j=1}^m$.

Рассмотрим теперь нормированные условия $\{{U}_j(y)=0\}_{j=1}^m$.
Для каждого $k=1,\ldots,m-1$ имеется не более двух краевых условий
порядка $k$. Будем обозначать $j(k)$ меньший из номеров краевых
условий порядка $k$, если условия такого порядка имеются.
Если для некоторого $k$ найдутся ровно два условия порядка $k$ с
номерами $j(k)$ и $j(k)+1$, то, заменяя ${U}_{j(k)}(y)=0$ и
${U}_{j(k)+1}(y)=0$ на их линейные комбинации, совершим
преобразования, приводящие к таким значениям коэфиициентов:
$\alpha_{j(k)}=\beta_{j(k)+1}=0$,
$\beta_{j(k)}=\alpha_{j(k)+1}=1$. Вычитая теперь из $U_n(y)=0$
для каждого $n<k(j)$ равенство
$\beta_nU_{j(k)}(y)+\alpha_nU_{j(k)+1}(y)=0$, добьёмся того,
чтобы $y^{(k)}(0)$ и $y^{(k)}(1)$ вcтречались только в условиях с
номерами $j(k)$ и $j(k)+1$. Если же для некоторого $k$ только
одно условие имеет порядок $k$, то, заменяя условия с номерами,
меньшими $j(k)$, на линейные комбинации их с условием
$U_{j(k)}(y)=0$, исключим  либо $y^{(k)}(1)$ (если
$\beta_{j(k)}\neq 0$), либо $y^{(k)}(0)$ (если $\beta_{j(k)}=0$)
из всех условий кроме $U_{j(k)}(y)=0$. После этих преобразований
краевые условия выглядят следующим образом:
$$
\begin{pmatrix}
\cdots& \cdot & 0 & \cdots & 0 & 0 & \cdots & \alpha_{k_1} & \beta_{k_1} & \cdots \\
       &\cdot   &  \cdot  &        &  \cdot & \cdot &&\cdot &\cdot &\\
&\cdot& 0&\cdots &0& 1&\cdots &0&0 &\\
&\cdot& 0&\cdots &1& 0&\cdots &0&0 &\\
&\cdot&\cdot&\cdot&\cdot&\cdot&\cdot&\cdot&\cdot &\\
\cdots & \alpha_k & \beta_k & \cdots &0 &0& \cdots & 0& 0 &\cdots\\
\cdots &\cdots &\cdots &\cdots &\cdots &\cdots &\cdots &\cdots
&\cdots &\cdots
\end{pmatrix}
\left( \begin{array}{c}
y(0)\\ y(1)\\ \cdot\\ \cdot \\ \cdot \\ y^{(m-1)}(0)  \\
y^{(m-1)}(1) \end{array} \right)=0
$$
Сделанные преобразования не изменили ни суммы порядков, ни
нормированность условий.

Теперь очевидно, что при совершении элементарных преобразований
строк в полученной матрице сумма порядков не уменьшится, а значит,
нормированные условия реализуют минимальную форму.
\end{proof}

Из этого утверждения очевидным образом получается следующее.
\begin{cons}
Набор порядков $k_1,\ldots,k_m$ нормированных условий не зависит
от того, какими преобразованиями проводилась нормировка.
\end{cons}

Теперь определим класс {\it регулярных} краевых условий.
Обозначим через $\omega_j, \; j=1,\ldots,m$ различные корни $m$-й
степени из $-1$ в таком  порядке, чтобы выполнялись неравенства
$$
\Re (\omega_1 e^{\frac{i\pi}{2m}}) < \Re (\omega_2
e^{\frac{i\pi}{2m}}) < \ldots < \Re (\omega_m
e^{\frac{i\pi}{2m}}).
$$
Определение даётся отдельно для чётных и нечётных $m$.

\begin{itemize}
\item[а)] $m$ нечётно; $m=2\mu-1$.

Нормированные краевые условия \eqref{normbc} называются
регулярными, если оба числа $\theta_0$ и $\theta_1$, определённые
равенством
$$\theta_0+\theta_1 s =
\begin{vmatrix} \alpha_1 \omega_1 ^{k_1}& \ldots &\alpha_1
\omega_{\mu-1} ^{k_1}& (\alpha_1+s\beta_1)\omega_{\mu} ^{k_1}&
\beta_1 \omega_{\mu+1} ^{k_1}&\ldots& \beta_1 \omega_{m} ^{k_1}\\
\alpha_2 \omega_1 ^{k_2}& \ldots &\alpha_2 \omega_{\mu-1} ^{k_2}&
(\alpha_2+s\beta_2)\omega_{\mu} ^{k_2}& \beta_2 \omega_{\mu+1}
^{k_2}&\ldots& \beta_2 \omega_{m} ^{k_2}\\
\cdot&\cdot&\cdot&\cdot&\cdot&\cdot&\cdot\\\alpha_m \omega_1
^{k_m}& \ldots &\alpha_m \omega_{\mu-1} ^{k_m}&
(\alpha_m+s\beta_m)\omega_{\mu} ^{k_m}& \beta_m
\omega_{\mu+1}^{k_m}&\ldots& \beta_m \omega_{m}
^{k_m}\end{vmatrix},$$ отличны от нуля.

\item[б)] $m$ чётно; $m=2\mu$.

Нормированные краевые условия \eqref{normbc} называются
регулярными, если по крайней мере одно из чисел $\theta_{-1}$ и
$\theta_1$, определённых равенством
$$\frac{\theta_{-1}}{s}+\theta_0+\theta_1 s =$$
$$=\begin{vmatrix} \alpha_1 \omega_1 ^{k_1}& \ldots &\alpha_1
\omega_{\mu-1} ^{k_1}& (\alpha_1+s\beta_1)\omega_{\mu} ^{k_1}&
(\alpha_1+\frac{1}{s}\beta_1)\omega_{\mu+1} ^{k_1}&
\beta_1 \omega_{\mu+2} ^{k_1}&\ldots& \beta_1 \omega_{m} ^{k_1}\\
\alpha_2 \omega_1 ^{k_2}& \ldots &\alpha_2 \omega_{\mu-1} ^{k_2}&
(\alpha_2+s\beta_2)\omega_{\mu} ^{k_2}&
(\alpha_2+\frac{1}{s}\beta_2)\omega_{\mu+1} ^{k_2}& \beta_2
\omega_{\mu+2}^{k_2}&\ldots& \beta_2 \omega_{m} ^{k_2}\\
\cdot&\cdot&\cdot&\cdot&\cdot&\cdot&\cdot&\cdot\\\alpha_m \omega_1
^{k_m}& \ldots &\alpha_m \omega_{\mu-1} ^{k_m}&
(\alpha_m+s\beta_m)\omega_{\mu} ^{k_m}&
(\alpha_m+\frac{1}{s}\beta_m)\omega_{\mu+1} ^{k_m} & \beta_m
\omega_{\mu+2}^{k_m}&\ldots& \beta_m \omega_{m}
^{k_m}\end{vmatrix}$$ отлично от нуля.
\end{itemize}

\section{Диссипативные краевые условия чётного порядка}
В этом параграфе мы докажем регулярность оператора чётного
порядка, порождённого $l_0(y)$ и диссипативными краевыми
условиями.

Будем писать $y_0^{(k)}$ вместо $y^{(k)}(0)$ и $y_1^{(k)}$ вместо
$y^{(k)}(1)$.

Центральное в этом параграфе утверждение звучит следующим образом.

\begin{thm}
Диссипативный оператор чётного порядка ($m=2n$) является
регулярным.
\end{thm}

\begin{proof}[Доказательство]
  Произведем нормировку краевых условий:
  \begin{equation}
  \label{normb}
  \mathit{U}_j(y)=b_0^j y_0^{(j)}+b_1^j y_1^{(j)}+\ldots =0,\;
  j=0,\ldots,2n-1
  \end{equation}
  $\mathit{U}_j(y)$ и $b_k^j$ --- векторы-столбцы высоты ${r_j \in
  \{0, 1, 2\}, \; \sum\limits_{j=0}^{2n-1}r_j=2n, \; k=1,2}$
  (многоточие здесь подразумевает производные младше $j$). Будем
  называть $y_k^{(j)}$ неизвестными. Имеем систему  линейных
  уравнений относительно этих неизвестных, которую можем решить,
  выразив главные неизвестные через свободные. Тогда краевые условия
  задают однозначное отображение из $\mathbb C^{2n}$ в
  $\mathbb C^{2n}$ --- каждому набору значений свободных
  неизвестных соответствует единственный набор главных так, что функция,
  значения которой в $0$ и $1$ вместе со значениями её производных
  до \mbox{$(2n-1)$--й} включительно совпадают с набором значений неизвестных,
  удовлетворяет краевым условиям.

  \begin{lemma} \label{lm1}

  В каждой паре $(y_k^{(j)},\, y_k^{(2n-1-j)})$, $k=0,1$ обе неизвестные не
  могут быть одновременно свободными, если условия диссипативны и
  $r_j\,\cdot \,r_{2n-1-j}\neq 1$.
  \end{lemma}

  \begin{proof}[Доказательство] Предполагаем противное. Пусть обе
  неизвестные в паре $(y_0^{(j)},\, y_0^{(2n-1-j)})$ свободные
  при некотором $j$. Запишем форму ${\Im{(\mathit{L}_0y,y)}}$,
  подставив вместо главных неизвестных   их выражения через
  свободные:
  \begin{multline*}
      \Im{(\mathit{L}_0y,y)}=(-1)^{n+j+1}\cdot
      \Im{(y_0^{(j)}\overline{y_0^{(2n-j-1)}})}+\eta
      |y_0^{(\operatorname{min}\{j,2n-j-1\})}|^2+\\
      +\alpha y_0^{(2n-1-j)}+\beta \overline{y_0^{(2n-1-j)}}+
      \gamma y_0^{(j)}+ \delta
      \overline{y_0^{(j)}} + \varepsilon,
  \end{multline*}
  где $\alpha\ldots\varepsilon$
  --- линейные комбинации свободных переменных за исключением
  рассматриваемой пары; $\eta$---число.

  Рассмотрим функции ${y \in D(\mathit{L}_0)}$, у которых значения
  производных, соответствующих свободным переменным, за исключением
  $y_0^{(j)}$ и $y_0^{(2n-1-j)}$ равны нулю. Форма
  ${\Im{(\mathit{L}_0y,y)}}$ на таких функциях принимает значения
  $$
  \Im{(\mathit{L}_0y,y)}=(-1)^{n+j+1}\cdot\Im{(y_0^{(j)}\overline{y_0^{(2n-j-1)}})}+\eta
      |y_0^{(\operatorname{min}\{j,2n-j-1\})}|^2
  $$
  и не является знакоопределенной. Противоречие. Лемма доказана.

  \end{proof}

  \begin{lemma} \label{lm2}
  Пусть краевые условия диссипативны, тогда

  \begin{equation}
  \label{rank_sum} r_j+r_{2n-1-j}=2 \quad  \forall j \in \{0,\ldots
  ,2n-1 \}.
  \end{equation}

  \end{lemma}

  \begin{proof}[Доказательство] Действительно, если ${r_j+r_{2n-1-j}<2}$ для
  некоторого $j$ (${r_j=0}$, ${r_{2n-1-j}=0 \text{ или } 1}$), то, очевидно,
  есть пара свободных неизвестных ${(y_k^{(j)},\, y_k^{(2n-1-j)})}$
  при \mbox{$k=0$ или $1$.} И в силу леммы \ref{lm1}, для диссипативных
  краевых условий не может выполняться ${r_j+r_{2n-1-j}<2}$ ни при каком $j$.

  Если же ${r_j+r_{2n-1-j}>2}$ при некотором $j$ и $\sum\limits_{l=0}^{2n-1}r_l=2n$,
  то найдётся такое $s$, что ${r_s+r_{2n-1-s}<2}$, а этого быть не может в
  силу сказанного выше. Лемма \ref{lm2} доказана.
  \end{proof}

  \begin{lemma} \label{lm3}
  Пусть заданы диссипативные краевые условия и
  нашлось такое $j$, что ${r_j=r_{2n-1-j}}$, тогда
  \begin{equation}
  \label{rank_1}
  b_0^{j}\overline{b_0^{2n-1-j}}=b_1^{j}\overline{b_1^{2n-1-j}}.
  \end{equation}

  \end{lemma}

  \begin{proof}[Доказательство] Сперва предположим, что среди этой четверки
  чисел нет нулей. Тогда рассмотрим неизвестные $y_1^{(j)}$ и
  $y_1^{(2n-1-j)}$ как свободные. Как и в доказательстве леммы \ref{lm1} положим
  значения всех свободных  неизвестных за исключением этой пары равными
  нулю. Тогда
  $$
     \Im{(\mathit{L}y,y)}=(-1)^{n+j}\cdot\Im{\left(y_1^{(j)}\overline{y_1^{(2n-j-1)}}
     \left[1-\frac{b_0^j}{b_1^j} \overline{\left(
     \frac{b_0^{(2n-j-1)}}{b_1^{(2n-j-1)}} \right)}\right]\right)}
     +\mu |y_1^{(\operatorname{min}\{j,2n-j-1\})}|^2,
  $$

  где $\mu$---некоторое число. Отсюда получаем требуемое.

  Рассмотрим оставшийся вариант: пусть $b_0^j=0, \, b_1^j\neq 0$.
  Для знакоопределенности $\Im{(\mathit{L}_0y,y)}$ на функциях
  $y$, у которых значения всех производных, соответствующих свободным
  неизвестным, отличным от выбранной четверки, нулевые, необходимо
  выполнение условий $b_1^{2n-1-j}=0, \, b_0^{2n-1-j}\neq 0$.
  \end{proof}

\begin{lemma} \label{lm5}

Пусть даны диссипативные нормированные краевые условия
$$
U_{j}(y)=b_0^j y^{(k_j)}(0)+b_1^j y^{(k_j)}(1)+\ldots=0, \;
j=1,\ldots,2n,
$$
тогда краевые условия, получаемые отбрасыванием младших
производных
$$
U^{(sa)}_{j}(y)=b_0^j y^{(k_j)}(0)+b_1^j y^{(k_j)}(1)=0, \;
j=1,\ldots,2n,
$$
являются самосопряжёнными.
\end{lemma}
\begin{proof} Это простое следствие лемм \ref{lm1}, \ref{lm2},
\ref{lm3}
\end{proof}
Поскольку в определении регулярности краевых условий участвуют
только коэффициенты при значениях старших производных, то
диссипативные условия $\{U_j(y)=0\}_{j=1}^{2n}$ и полученные по
ним самосопряжённые $\{U^{(sa)}_{j}(y)=0\}_{j=1}^{2n}$ регулярны
одновременно. В этом месте сошлёмся на работу S.Salaff
\cite{Sal}, в которой установлена регулярность самосопряжённых
дифференциальных операторов чётного порядка. Доказательство
теоремы закончено.
\end{proof}
\section{Диссипативные краевые условия нечётного порядка.}

Диссипативные операторы могут не являться регулярными. Для
доказательства этого приведём следующий пример.

\begin{example}

$$
l(y)=(-i)^{2n-1} y^{(2n-1)}=(-1)^n i y^{(2n-1)}
$$

Краевые условия возьмем такими:

$$
y_0^{(2n-2)}=y_1^{(2n-2)}=\ldots =y_0^{(n)}=y_1^{(n)}=y_1^{(n-1)}=0
$$

При заданных краевых условиях
$\Im{(\mathit{L}_0y,y)}=\frac{|y_0^{(n-1)}|^2}{2},$
\end{example}

Элементарная проверка показывает, что число $\theta_0$, участвующее в определении регулярности равно нулю. Поэтому регулярными такие краевые  условия не являются.

\vspace{5mm} Автор благодарит А.~А.~Шкаликова за постановку
задач, обсуждение работы и ценные замечания. В частности, ему принадлежит формулировка теорем 1 и 2.

\end{document}